\documentclass[reqno,11pt]{amsart}

\usepackage{mathrsfs,amssymb, latexsym,amsmath,mytitle}

\newtheorem{thm}{Theorem}[section]

\newtheorem{lemma}[thm]{Lemma}

\theoremstyle{remark}
\newtheorem{rmk}{Remark}

\newcommand{\N}{\mathrm{N}}
\newcommand{\End}{\mathrm{End}}
\newcommand{\tr}{\mathrm{Tr}}
\newcommand{\nd}{\mathrm{and}}

\newcommand{\fo}{\mathrm{for}}
\newcommand{\sn}{\mathrm{sgn}}

\newcommand{\fg}{\mathfrak{g}}
\newcommand{\fsl}{\mathfrak{s}\mathfrak{l}}

\newcommand{\VF}[2]{{\mathcal V}_{{#1}, {#2}}}
\newcommand{\VFR}[2]{{\mathcal V}_{\! {#1}, {#2}}^{\scriptscriptstyle Y}}
\newcommand{\VFRp}[3]{{\mathcal V}_{\! {#1}, {#2}}^{{\scriptscriptstyle Y}#3}}

\newcommand{\fh}{\mathfrak{h}}
\newcommand{\bbz}{\mathbb{Z}}

\newcommand{\bbr}{\mathbb{R}}
\newcommand{\bbc}{\mathbb{C}}
\newcommand{\sS}{\mathscr{S}}

\newcommand{\slz}{\mathrm{SL}(2,\bbz)}
\newcommand{\zo}{\bbz[\omega]}
\newcommand{\hs}{\hspace}
\newcommand{\ur}{U^{(r)}}
\newcommand{\sltz}{\mathrm{SL}(2,\bbz)}
\newcommand{\psu}[1]{\mathrm{PSU}(#1)}

\newcommand{\old}[1]{}

\newcommand{\subjectclass}[1]{
\renewcommand{\thefootnote}{} 
\footnote{2000\ \textit{Mathematics Subject Classification.} \textrm{#1}}
\addtocounter{footnote}{-1}
\renewcommand{\thefootnote}{\arabic{footnote}} }

\address{Department of Mathematics, The Ohio State University, Columbus, OH 43210, USA}
\email{qichen@math.ohio-state.edu} 
\address{Department of Mathematics, The Ohio State University, Columbus, OH 43210, USA}
\email{kerler@math.ohio-state.edu} 

\renewcommand{\thefootnote}{}

\begin{document}

\institute{}

\title{Higher Rank TQFT Representations of $\mathrm{SL}(2,\bbz)$ are Reducible}
\subtitle {Examples of Decompositions and Embeddings.}

\bigskip

\author{Qi Chen  \ \ and \ \  Thomas Kerler}

\maketitle

\begin{abstract}
In this article we give examples which show that the TQFT representations
of the mapping class groups derived from quantum $\psu{N}$ for $N>2$ are generically decomposable. 
One general decomposition of the representations 
is induced by the symmetry which exchanges
$\psu{N}$ representation labels by their conjugates. The respective summands of a given
parity are typically  still reducible into many further components.

\medskip

Specifically, we give an explicit basis for an  irreducible direct summand in the $\sltz$
representation obtained from quantum $\psu{3}$ when the order of the root of unity is a prime $r\equiv 2 \mod 3$.
We show that this summand is isomorphic to the respective  $\psu{2}$ representation.
\end{abstract}

\subjectclass{Primary: 57R56; Secondary: 57M27,11F27}


\section{Introduction}\label{intro}

\subsection{Background and Motivation}  

One of the most interesting features of a topological quantum field theory 
(TQFT)  in dimension three is that it naturally implies
projective representations of mapping class groups of surfaces.
The question which motivates this article is how the decomposition structure of these 
representations changes as one passes from TQFTs constructed from 
$\fsl_2$ to TQFTs constructed from higher rank Lie algebras.

More precisely,  a TQFT can be constructed from the representation theory 
of a finite dimensional quasi-triangular Hopf algebra. 
Particularly, for any complex simple Lie algebra $\fg$ and any root of 
unity $\zeta$  such a Hopf algebra is given by 
a quantum group  $U_\zeta(\fg)'\,$, and, hence, defines an associated TQFT which we denote by  $\;\VF{\fg}{\zeta}\,$.

Basis elements of the vector spaces of this TQFT are typically enumerated by graphs or links whose edges are
labeled by weights of $\,\fg\,$ in the associated Weyl alcove. Restricting $\,\VF{\fg}{\zeta}\,$
to subspaces in which basis elements are labeled by only roots of $\,\fg\,$, one naturally
obtains a  sub-TQFT denoted by $\,\VFR{\fg}{\zeta}\,$. It turns out that  $\,\VF{\fsl_N}{\zeta}\,$ differs
from  $\,\VFR{\fsl_N}{\zeta}\,$ only by tensoring with a trivial TQFT based on $U(1)$  \cite{blan00,blan05}.
Commonly,  $\;\VFR{\fsl_N}{\zeta}\,$ is also referred to as the TQFT associated to quantum $\psu{N}$.

In the case where  $\;\fg =\fsl_2\;$ and $\;\zeta\;$ is of prime order Roberts proved in \cite{roberts}
that the mapping class group representations obtained  from $\,\VFR{\fsl_2}{\zeta}\,$   are indeed irreducible.
Irreducibility is a crucial feature of the $\,\fsl_2$-TQFT
in a number of applications such as the 
generation of integral bases from cyclic vectors or mutation invariance with
respect to the central hyperelliptic involutions in genus $\leq 2\,$.
Yet, hardly anything is known about the reducibility question for  representations obtained from 
$\,\VFR{\fg}{\zeta}\,$ for higher rank Lie algebras $\,\fg\,$ or for   $\,\zeta\,$'s of non-prime orders.

The main result  of this article is  that irreducibility fails in several ways for higher rank Lie algebras. More precisely, we provide 
explicit decompositions of the representations of  the mapping class group of the torus $\Gamma_1\cong \sltz$ for $\,\fg=\fsl_3\,$
 but $\,\zeta\,$ still of prime order and identify respective summands with representations coming from $\fsl_2\,$. 
The genus one case reveals  already a fair amount of the structure of the TQFT on general cobordisms, 
suggesting reducibilities also for higher genera.

In fact, representations of the modular group  $\sltz$ have been a topic of intensive study in  conformal field theory
as well as analytic number theory. In particular, the closely related CFT modular invariants for $\,\fg=\fsl_3\,$ have
 been analyzed and classified, for example,  in \cite{bcir,gannon,itzyk,rtw}, and we will use symmetries already 
mentioned in \cite{itzyk,rtw}.

The computations of this paper also serve as explicit examples and a starting point 
for the more general approach taken in \cite{ck}. There we use methods and techniques 
previously developed for the study of theta functions in analytic number theory in order to 
give the complete decompositions of the $\sltz$ representations obtained from
$\,\VFR{{\mathfrak g}}{\zeta}\,$ into irreducible representation types of $\sltz$
for all Lie algebras  $\mathfrak g$ of rank two when  $\zeta$ is of order $p^{\lambda}$ 
with $p$ an odd prime. Further work will investigate in how far decompositions found 
here will extend to higher genus  representations.

\bigskip

\subsection{Statement of Results} Before turning to our main result for $\,\fg=\fsl_3\,$, we 
consider the  $\,{\mathbb Z}/2$-symmetry  obtained by replacing weights or color  
labels by their conjugates.  It  naturally endows a large family of  TQFTs constructed from
 higher rank Lie algebras with a respective $\,{\mathbb Z}/2$-grading.

\begin{thm}\label{main-thm-1}
Suppose  $\,N>2\,$ and  $\,\zeta\,$ is of  order $k> 2N$ , with $k$ and $N$ coprime. Then the TQFT $\VFR{\fsl_N}{\zeta}$ 
(as in \cite{y2,lickorish,blan00}) commutes with a   $\,{\mathbb Z}/2$-action on each $\,\VFR{\fsl_N}{\zeta}(\Sigma)\,$ for any closed, compact, orientable surface $\Sigma$. 
Consequently,  the TQFT is $\bbz/2$-graded with a decomposition as follows:
$$ 
\VFR{\fsl_N}{\zeta}\,=\,\VFRp{\fsl_N}{\zeta}{+}\,\oplus\,\VFRp{\fsl_N}{\zeta}{-}\, .
$$ 
The  $\,{\mathbb Z}/2$-action on  $\,\VFR{\fsl_N}{\zeta}(\Sigma)\,$ is non-trivial whenever  $\Sigma$ is not the union of spheres $S^2$.
\end{thm}

In particular, the representation of the mapping class group of $\Sigma$ induced by $\,\VFR{\fsl_N}{\zeta}\,$ is obviously
reducible   into the two graded parts.  The basic parity splitting from Theorem~\ref{main-thm-1}  is, however, far from the only 
decomposition of the mapping class group representations in the higher rank case. 
 
\medskip

In fact,  the summands $\,\VFRp{\fsl_N}{\zeta}{\epsilon}\,$ are generally still further reducible into many more components. 
Particularly, in the  genus one case for $\zeta$ of prime order only a few representation types of the modular group $\sltz$
can occur which, necessarily, will repeat themselves more and more often in decompositions for higher rank Lie algebras. 

The main theorem of this paper makes this phenomenon explicit in the concrete example of the $\sltz$-representation of
$\psu{3}$ type, which turns out to contain summands of representations of $\psu{2}$-type as well as other representation types.

\begin{thm}\label{main-thm-2}
Suppose that $r$ is an odd prime and $\zeta = \exp(\frac{2\pi \sqrt{-1}}r)$. Assume further  $r\equiv 2 \mod 3$ and let $\epsilon=(-1)^{\frac {r+1}6}$. 
Then the $\sltz$ representation obtained from $\,\VFR{\fsl_2}{\zeta^4}\,$ occurs 
as a summand in  the $\sltz$ representation obtained from $\,\VFRp{\fsl_3}{\zeta}{\epsilon}\,$.
\end{thm}

The proof will  identify a specific basis for the respective submodule of  the  $\sltz$ representation of $\psu{3}$ type.
 The equivalence is shown by explicitly computing that the restricted matrix coefficients of the  $\sltz$ generators acting on this
subspace  coincide with the respective coefficients for the $\psu{2}$ type representation of $\sltz$. The decomposition is derived
from additional symmetries of the root lattice, given that the latter can be viewed as a quadratic extension over $\mathbb Z$ 
with a dihedral automorphism group.  

\medskip
\medskip

{\em Acknowledgment.} The first author would like to thank Tom Cusick for pointing out the reference \cite{jacobsthal} to him.
He also would like to thank Gregor Masbaum for helpful discussions.

\bigskip


\section{Proof of Theorem \ref{main-thm-1}}\label{skein}

\begin{proof}
The  $\bbz/2$-action asserted in Theorem~\ref{main-thm-1} can be constructed for all TQFTs, ${\mathcal V}_{\mathscr S}$,  for which the vector space
${\mathcal V}_{\mathscr S}(\Sigma)$ associated to surface $\Sigma$ is obtained from the skein module  $\sS(H)$  of a handle body, $H$,
with $\;\Sigma=\partial H\,$ and for which the following further assumptions hold. 

Firstly, 
we assume  $\sS(S^3)={\mathbb C}\,$ so that for any link $Q\subset S^3$ an invariant $P(Q)\in{\mathbb C}$ is
defined by  $\{Q\}=P(Q)\{\emptyset\}$, where $\{Q\}\in \sS(S^3)$ denotes the skein class of the link.  We also assume  that the relations underlying the skein module are orientation invariant so that
there is a well defined  $\bbz/2$-action on   $\sS(H)$ obtained by mapping the class $\{L\}\in \sS(H)$ 
of  each link $L$ to the class $\{\overline L \}$ of the 
link $\overline L$ with reversed orientation. Note, that this implies $P(Q)=P(\overline Q)$ for the invariant since $\overline{\emptyset}=\emptyset$. 
Finally, we assume that the  special element $\Omega$ in $\sS(S_1\times D^2)$  implementing 0-surgery along the core of the full torus is also
invariant under orientation reversal. 

Recall, that the vector space ${\mathcal V}_{\mathscr S}(\Sigma)$ is obtained as a subquotient of  
$\sS(H)$ as follows.  Consider an unknotted embedding of $H$ into $S^3$ such that its complement
is again a handle body $H'$. Given that $S^3=H\cup_{\Sigma}H'$ a link $L\subset H$ and 
a link $L'\subset H'$ obviously combine to a link $L\cup L'\subset S^3$, so that we obtain 
the following bilinear form on the level of skein classes $\{L\}$ of links:
$$
\langle\ , \ \rangle : \sS(H)\times\sS(H')\longrightarrow \bbc\,:\hs{3ex} (\{L\}, \{L'\})\,\mapsto \, \langle \{L\}, \{L'\}\rangle := 
P(L\cup L').
$$ 
The space  ${\mathcal V}_{\mathscr S}(\Sigma)$  is now defined by dividing  $\sS(H)$ by the null-space of this form. 
Since   $\langle \{L\}, \{L'\}\rangle= \langle \{\overline L\}, \{\overline{L'}\}\rangle$ the null-space is preserved under the 
$\bbz/2$-action so that it is well defined on  ${\mathcal V}_{\mathscr S}(\Sigma)\,$. 
For a framed link $L$ in $H$, we denote by $[L]$   the corresponding element in $\,{\mathcal V}_{\mathscr S}(\Sigma) \,$,
and by $\,\varphi\in\mathrm{End}({\mathcal V}_{\mathscr S}(\Sigma))\,$ the involution given by $\varphi([L])=[\overline L]$
which generates the $\bbz/2$-action.
For a disconnected surface $\Sigma=\bigsqcup_{j=1}^m\Sigma_j$ the $\bbz/2$-action on 
$\,{\mathcal V}_{\mathscr S}(\Sigma)\,=\,\bigotimes_{j=1}^m{\mathcal V}_{\mathscr S}(\Sigma_j)\,$ is given by $\varphi^{\otimes m}\,$. 

In order for ${\mathcal V}_{\mathscr S}(\Sigma)$  to be $\bbz/2$-graded the TQFT needs to be equivariant with respect to this 
action. That is, $\varphi$ has to commute with ${\mathcal V}_{\mathscr S}(M): {\mathcal V}_{\mathscr S}(\Sigma_{S})\to{\mathcal V}_{\mathscr S}(\Sigma_T)\,$ 
for any cobordism $\,M:\Sigma_S\to\Sigma_T\,$. It suffices to show this for connected $M$. Suppose $H_T$ and $H_S$ are collections of handle bodies
with $\Sigma_{S(T)}=\partial H_{S(T)}$, $H'_{S(T)}$ are respective complements in $S^3$, and $H_T^{\#}$ is the interior connected sum of components of $H_T\,$.
It is well known that $M$ can be obtained by surgery on a framed link $\mathscr{L}\subset N:=H_T^{\#}-H_S\,$ (see, e.g., \cite{k99} for details and link calculus). 

The spaces  ${\mathcal V}_{\mathscr S}(\Sigma_{S})$ and  ${\mathcal V}_{\mathscr S}(\Sigma_T)^*$ are spanned by classes of links $L_S\in H_S$ and 
$L_T\in H_T'$.  A general matrix element is thus given by 
\begin{equation}\label{eq-VMme}
\langle [L_T], {\mathcal V}_{\mathscr S}(M). [L_S]\rangle\;=\;P(L_T\cup \mathscr{L}(\Omega)\cup L_S)\;
\end{equation}
where we use that $H_T'\cup N\cup H_S=S^3$ and where $\mathscr{L}(\Omega)$ denotes the skein in which each component of $\mathscr{L}$ is 
substituted by $\Omega$. 

It is now obvious from   (\ref{eq-VMme})  that a matrix element of $\varphi\circ  {\mathcal V}_{\mathscr S}(M)$   differs from the one for 
$ {\mathcal V}_{\mathscr S}(M)\circ \varphi$  by orientation reversal on $\mathscr{L}(\Omega)$. However, since we assumed that
$\Omega$ is invariant under orientation reversal equality follows. 

\medskip

It remains to show that the TQFT as defined in  is  \cite{y2,lickorish,blan00} based on the HOMFLYPT polynomial  
fulfills the assumptions made at the beginning of the proof. Note, that the modularity condition needed to construct 
a TQFT from the skein theory was proved in  \cite{lick98} only for both $k$ and $N$ prime (see Lemma 4.2), but was later shown
in  \cite{blan00} to apply also to the more general situation with $(k,N)=1$ (see Lemma 2.9 and Theorem 2.11).

The fact that the HOMFLYPT skein module of $S^3$  is one dimensional is what makes
the HOMFLYPT polynomial of links well defined. It readily follows from inspection of the HOMFLYPT skein relation
diagrams that the respective skein ideal  is invariant under orientation reversal. Specifically, orientation reversal
yields the same diagrams as a rotation by $\pi$ in the plane of projection for the standard presentation of the relations. 

For the verification of the  last assumption, recall that the element $\Omega$ is (proportional to) a
 sum $\sum_{\lambda\in\Gamma}\Delta_{\lambda}\mathbf{e}_{\lambda}$, where $\Gamma$ is a set
of Young diagrams with $N-1$ rows and at most $k$ boxes in each row, 
$\mathbf{e}_{\lambda}$ the torus skein on $|\lambda|$ strands with inserted
Wenzl/Yokota idempotent, and $\Delta_{\lambda}$ the associated quantum dimension given by a standard 
embedding of $\mathbf{e}_{\lambda}$ 
in $S^3$.  Rigidity of the underlying  Hecke category and Proposition 2.6 of  \cite{blan00} imply that $\overline{\mathbf{e}_{\lambda} }
=\mathbf{e}_{\lambda^*}$, where $\lambda^*$($=\lambda_1^N/\lambda$) is the dual diagram 
so that also
  $\Delta_{\lambda}=\Delta_{\lambda^*}$. Since $\Gamma^*=\Gamma$ this implies $\overline\Omega=\Omega$ as assumed.

Finally, in order to prove that the $\bbz/2$-action is non-trivial on $\,\VFR{\fsl_N}{\zeta}(\Sigma)\,$ for $\,\Sigma\neq S^2\,$ consider
the skeins given by ${\mathbf{e}_{\lambda} }$ and $\mathbf{e}_{\lambda^*}$ along the core of one of the handles of $H$ with
$\Sigma=\partial H$. These yield vectors $v_{\lambda}$ and  $v_{\lambda^*}$ in $\,\VFR{\fsl_N}{\zeta}(\Sigma)\,$ which are mapped to
each other under the  $\bbz/2$-action. Non-triviality of the latter thus follows if we show that these vectors are not colinear 
for $\lambda\neq\lambda^*\,$. 

To this end construct a skein $\mathbf{f}_{\lambda}$ in the complementary handle body $H'$ as follows. Push the meridian of
the handle chosen for  ${\mathbf{e}_{\lambda} }$  off into $H'\subset S^3$ and replace it by $\Omega$. Push off a meridian of 
$\Omega$ in $H'$ and replace it with ${\mathbf{e}_{\lambda} }$ with orientation opposite  to that of the skein in $H$. Denote the respective covector $w_{\lambda}$  in $\,\VFR{\fsl_N}{\zeta}(\Sigma)^*\,$.
It follows from modularity that for $\lambda\neq\lambda^*\,$ we have   $w_{\lambda}\circ v_{\lambda^*}=0$ but 
$w_{\lambda}\circ v_{\lambda}\neq 0$ for $\lambda\in\Gamma$. Hence   $v_{\lambda}$ and  $v_{\lambda^*}$ are not colinear and 
the $\bbz/2$-action is non-trivial.
\end{proof}

We expect Theorem \ref{main-thm-1} to apply also to TQFTs based on Lie algebras $\mathfrak{so}_{4n+2}$ and $E_6$ 
as well as other modular categories with non-trivial duals. A proof would require either
finding a skein theoretic formulation that fulfills the assumptions in the above proof, 
or proving respective symmetries of structure constants (such as 6-j symbols) of a
corresponding colored graph calculus.

\section{
Modular  $\fsl_3$ Root Lattice, Unfolded $\sltz$ Representations, and Their Symmetries.}\label{unf}
%
The automorphism group of the root lattice of   $\fsl_3$  modulo a prime $r$ with $r\equiv 2\mod 3$ is a
dihedral group of order $2(r+1)$, containing the ordinary Weyl group $W\cong S_3$. In this section we explicitly
describe the symmetries and decompositions this automorphism group entails for the associated unfolded (Weil) representation of 
$\sltz$ acting on the space of functions on the $r$-modular root lattice of  $\fsl_3$.

\subsection{\boldmath The  root lattice of $\fsl_3$}\label{root}
We fix a Cartan subalgebra $\fh$ of $\fsl_3$ and basis roots $\alpha_1$ and $\alpha_2$ in its dual space $\fh^*$. Let
$\fh^*_\bbr$ be the $\bbr$-subspace of $\fh$ spanned by $\alpha_1$ and $\alpha_2$. 
The root lattice $Y$ of $\fsl_3$ is the $\bbz$-span of $\alpha_1$ and $\alpha_2$. 
Define an inner product $(\cdot |\cdot)$ on
$\fh_\bbr^*$ such that $(\alpha_i|\alpha_j)=a_{ij}$ where $(a_{ij})$ is the Cartan matrix of $\fsl_3$.
The fundamental weights $\lambda_1$ and $\lambda_2$ are dual to the basis roots respect to this inner product, i.e.
$(\lambda_i|\alpha_j) = \delta_{ij}$. The weight lattice $X$ of $\fsl_3$ is the $\bbz$-span of $\lambda_1$ and $\lambda_2$.

Let $\omega$ be the third root of unity $(-1+\sqrt{-3})/2$. Then we can identify $X$ with the ring $\zo$ such that
$\lambda_1$ and $\lambda_2$ are identified with 1 and $-\omega^2$ respectively, see \cite{rtw}. 
The root lattice $Y$ is $(1-\omega)\zo$ under this identification. Since we
will only consider the root lattice, we may identify $Y$ with $\zo$ directly via the maps
$$
Y\to \zo \hs{3ex} \mathrm{with} \hs{3ex} \alpha_1 \mapsto 1 \hs{3ex} \textrm{and} \hs{3ex} \alpha_2 \mapsto \omega.
$$
Under this identification, $(z|z') = \tr(\bar z z')$ and $|z|^2 = 2\N(z)$ 
where $\bar z$ is the complex conjugation of $z$. 
The norm and trace $\zo\to \bbz$ are defined by
$$
\N(z) = z \bar z \hs{5ex} \mathrm{and} \hs{5ex} \tr(z) = z + \bar z.
$$
The Weyl group $W$ ($\cong S_3$) action on
$Y$ induces action on $\zo$:
$$
z\stackrel{w}{\rightarrow} z, \omega z, \omega^2 z, -\bar z, -\omega\bar z, -\omega^2\bar z
$$
where the first three are the orbit of even permutations.

\subsection{\boldmath Unfolded representations of $\slz$}\label{unfolded}
We denote by $S$ and $T$ the usual generators of  $\slz$ given by the matrices 
$$
S = \left(\begin{array}{cc}
0 & -1\\
1 & 0
\end{array}
\right) \hs{5ex} \mathrm{and} \hs{5ex} 
T = \left(\begin{array}{cc}
1 & 1\\
0 & 1
\end{array}
\right)\;.
$$
They obey the  relations $(ST)^3 = S^2$ and $S^4 = I$. 

We will  further assume for the remainder of this article that    $r$ is an odd prime   equal to 2 mod 3 so that
$r=6n-1$ for some $n$.  
The $r$-modular root lattice $Y/r$ is identified with  $R_r = \zo/r$, which is a field for   the assumed case $r\equiv 2 \mod 3$ (but {\em not} for $r$ is 1 mod 3).
For any $z\in\zo$ we write $(z \mod r) \in R_r$ simply as $z$.
The norm and trace on $\zo$  descend to norms and traces  $R_r\to \bbz/r\,$,  for which we will use the same notation.

An unfolded representation of $\sltz$ is defined with respect to an $r$-th root of unity $\zeta$.  We will restrict ourselves to the case where $\zeta$ 
is one of the roots closest to unit, that is, 
$$
\zeta = \exp(\frac{2\pi\sqrt{-1}}r)\;.
$$
Let $\hat V_r$ be the $\bbc$-vector space with basis $\hat B_r = \{e_z : z\in R_r\}$. We define a projective representation as follows:
\begin{equation}\label{phir}
\hat\phi_\zeta : \slz\to \End(\hat V_r) \hs{3ex} \textrm{with} \hs{3ex} \hat\phi_\zeta(S) = \hat S_\zeta \hs{3ex} \textrm{and}
\hs{3ex} \hat\phi_\zeta(T) = \hat T_\zeta,
\end{equation}
with matrix elements give by  
$$
\;(\hat S_\zeta)_{xy} = \zeta^{\tr(\bar xy)}\;\qquad\mbox{and}\qquad\;(\hat T_\zeta)_{xy} = \delta_{xy}\zeta^{-\N(x)}\;.
$$
\subsection{\boldmath The first symmetry of ${\hat\phi_\zeta}$}\label{1st_symmetry}
The field $R_r$ has $r^2$ elements. Its nonzero elements form a cyclic group $R_r^*$ of order $r^2-1$ under multiplication.
It contains the subgroup 
$$\ur = \{a\in R_r^* : \N(a)\equiv 1 \mod r\}$$
which is of order $6n = r+1$  \cite{ir} . We fix a generator $u$ of $\ur$ such that
$u^{2n}= \omega$ and $u^{3n} = -1$.

We define a $\ur$-action on $\pm \hat B_r$ by $u(e_z) = - e_{uz}$. It induces a $\ur$-action on $\hat V_r$ by linearity. It is
easy to check that this $\ur$-action on $\hat V_r$ commutes with the $\slz$-action
defined in Equation~(\ref{phir}). Therefore, the orbit sums of $\pm \hat B_r$ span an invariant subspace $\bar V_1$ of $\hat\phi_\zeta$. We follow \cite{rtw} to fix a basis for $\bar V_1$.
Let $\rho$ be an element in $R_r^*$ such that $\N(\rho)=\rho\bar\rho$ is not a square in $\bbz/r$.
Then 
$$
R_r^* = \{a\rho^bu^j : 1\le a\le \frac{r-1}{2}, b \in\{ 0, 1\}, 0\le j\le r\}.
$$
Recall that $r\equiv 2 (3)$. When $r\equiv 1 (4)$ we will take $\rho = \omega - \bar\omega$ since in this case $\N(\rho) = -3$,
which is not a square mod $r$. When $r\equiv 3 (4)$ we will take $\N(\rho)=-1$. We may further require
that $\bar\rho = u\rho$ when $r\equiv 3 (4)$.
 
For $1\le a\le \frac{r-1}{2}$, let
$$
e^+_a=\sum_{j=0}^r u^j(e_a) \hs{5ex}\nd\hs{5ex}
e^-_a =\sum_{j=0}^r u^j(e_{a\rho}).
$$
Obviously $\{e^+_a, e^-_a\}$ is a basis of $\bar V_1$.

\subsection{\boldmath The second symmetry of $\hat\phi_\zeta$}\label{2nd_symmetry}
Define an action of the  Weyl  $\,W$  on the set $\,\pm\hat B_r\,$ by $\;w(e_z)=\sn(w)e_{w(z)}$.
By linearity  it  induces a $W$-action on $\hat V_r$. It is easy to check that this $W$-action on $\hat V_r$ commutes with the $\slz$-action
defined in Equation~(\ref{phir}). Therefore, the orbit sums of $\pm \hat B_r$ span an invariant subspace $\bar V_2$ of the representation
$\hat\phi_\zeta$. We will follow \cite{rtw} again to fix a basis for $\bar V_2$. To do this we need to consider the cases $r\equiv 1 (4)$ and $r\equiv 3 (4)$ separately:

\subsubsection{$\mathbf{r\equiv 1 (4)}$}\label{odd}
If $r\equiv 1 (4)$ then $n=(r+1)/6$ is odd. 
For any $1\le a\le \frac{r-1}{2}$, let
\begin{eqnarray*}
e^{+o}_{a,j} = \sum_{w\in W}w(e_{au^{2j}}) & \fo & \hs{1ex}1\le j\le n, \\
e^{-o}_{a,j} = \sum_{w\in W}w(e_{a\rho u^j}) & \fo & \hs{1ex} 1\le j\le n-1.
\end{eqnarray*}
Then according to Equations~(3.22) and (3.23) in \cite{rtw} (where $n$ is denoted $q$),
$\{e^{+o}_{a,j}, e^{-o}_{a,j}\}$ is a basis of $\bar V_2$. Note
that $e^{-o}_{a,n}=0$.

\begin{lemma}\label{inv_odd}
Suppose $r\equiv 1 (4)$.
Let $V_r$ be the $\bbc$-vector subspace of $\hat V_r$ with basis $B_r = \{e^+_a: 1\le a\le \frac{r-1}{2}\}$.
Then $V_r$ is an $\slz$-invariant space. 
\end{lemma}
\begin{proof}
For any $1\le a\le \frac{r-1}{2}$,
\begin{eqnarray*}
e^+_a &=& \sum_{j=0}^{n-1}(-1)^j(e_{au^{2j}} + e_{au^{2n+2j}}+e_{au^{4n+2j}}-e_{au^{n-2j}}-e_{au^{3n-2j}}-e_{au^{5n-2j}}) \\
&=& \sum_{j=0}^{n-1}(-1)^j\sum_{w\in W}w(e_{au^{2j}}) = \sum_{j=0}^{n-1}(-1)^j e^{+o}_{a,j}.
\end{eqnarray*}
Therefore, $V_r\subset \bar V_1\cap \bar V_2$.
$$
\hat S_\zeta(e^+_a) = \sum_{b=1}^{(r-1)/2}c^+_be^+_b + \sum_{b=1}^{(r-1)/2}c^-_be^-_b
= \hat S_\zeta(\sum_{j=0}^{n-1}(-1)^j e^{+o}_{a,j})
$$
where $c^+_b$ and $c^-_b$ are complex numbers. Since the last sum does not contain $e_{b\rho u^n}$ terms,
$c_b^- = 0$. It remains to notice that $\hat T_\zeta(e^+_a) = \zeta^{-a^2} e_a^+$.
\end{proof}

Denote the restriction of $\hat\phi_\zeta$ to $V_r$ by $\phi_\zeta$. Also denote $\phi_\zeta(S)$ by $S_\zeta$ and
$\phi_\zeta(T)$ by $T_\zeta$. It is easy to see that in the basis $B_r$, 
\begin{equation}\label{eq:PSU(3)}
(S_\zeta)_{ab} = \sum_{j=0}^r(-1)^j\zeta^{ab\tr(u^j)} \hs{3ex}\nd\hs{3ex} (T_\zeta)_{ab}=\delta_{ab}\zeta^{-a^2}.
\end{equation}
Obviously, $(S_\zeta)_{ab} = (S_\zeta)_{ba} = -\overline{(S_\zeta)}_{ab}$.

\subsubsection{$\mathbf {r\equiv 3 (4)}$}\label{even}
If $r\equiv 3 (4)$ then $n$ is even. 
For any $1\le a\le \frac{r-1}{2}$, let
\begin{eqnarray*}
e^{+e}_{a,j} = \sum_{w\in W}w(e_{au^{\frac{n}{2}+j}}) & \fo & \hs{1ex}1\le j\le n-1, \\
e^{-e}_{a,j} = \sum_{w\in W}w(e_{a\rho u^{2j}}) & \fo & \hs{1ex} 1\le j\le n.
\end{eqnarray*}
Then according to Equations~(3.26) and (3.27) in \cite{rtw}, $\{e^{+e}_{a,j}, e^{-e}_{a,j}\}$ is a basis of $\bar V_2$. Note
that $e^{+e}_{a,n}=0$.

\begin{lemma}\label{inv_even}
Suppose $r\equiv 3 (4)$.
Let $V_r$ be the $\bbc$-vector subspace of $\hat V_r$ with basis $\{e^-_a: 1\le a\le \frac{r-1}{2}\}$.
Then $V_r$ is an $\slz$-invariant space. 
\end{lemma}
The proof is very similar to that of Lemma~\ref{inv_odd} hence is omitted. With the similar notation as in Section~\ref{odd}, we have
$$
(S_\zeta)_{ab} = \sum_{j=0}^r(-1)^j\zeta^{-ab\tr(u^j)} \hs{3ex}\nd\hs{3ex} (T_\zeta)_{ab}=\delta_{ab}\zeta^{a^2}.
$$
Note that we take $\N(\rho)=-1$ and $\bar \rho = u\rho$.
Obviously $(S_\zeta)_{ab} = (S_\zeta)_{ba} = \overline{(S_\zeta)}_{ab}$.

\begin{rmk}\label{rk}
Note that the formulas of $(S_\zeta)_{ab}$ and $(T_\zeta)_{ab}$ when $r\equiv 3 (4)$ are conjugate to those when
$r\equiv 1 (4)$.
\end{rmk}
\section{Quantum $\psu{3}$ Representations of $\slz$ and Proof of Theorem~\ref{main-thm-2} }\label{psun}
%
In this section we identify the summands of the unfolded representations from Section~\ref{unf}  with $\sltz$ representations 
arising from  $\psu{3}$ and   $\psu{2}$ TQFTs by explicit computations and comparisons of matrix elements. 
\subsection{\boldmath Quantum $\psu{2}$ representations}\label{PSU(2)}
For any odd prime $r$, let $V_r'$ be the $\bbc$-vector space with basis $B'_r=\{e'_a: 1\le a\le \frac{r-1}{2}\}$.
Let the (projective) quantum $\psu{2}$ representation of $\slz$ be
$$
\phi'_\xi : \slz\to \End(V'_r) \hs{3ex} \textrm{with} \hs{3ex} \phi'_\xi(S) = S'_\xi \hs{3ex} \textrm{and}
\hs{3ex} \phi'_\xi(T) = T'_\xi.
$$
Recall from \cite{le2} that for a primitive $r$th root of unity $\xi$,
$$
(S'_\xi)_{ab} = \xi^{2ab-a-b} - \xi^{-2ab+a+b-1} \hs{5ex}\nd\hs{5ex} (T'_\xi)_{ab}=\delta_{ab}\xi^{-a(a-1)}
$$
in the basis $B'_r$. Since $(4,r)=1$ we can take $\xi = \zeta^4$. Then
$$
(S'_\xi)_{ab} = \zeta^{2(2a-1)(2b-1)} - \zeta^{2(2a-1)(1-2b)} \hs{5ex}\nd\hs{5ex} (T'_\xi)_{ab}=\delta_{ab}\zeta^{-(2a-1)^2}.
$$
After reordering basis elements and possibly multiplying $-1$ we may assume
\begin{equation}\label{eq:PSU(2)}
(S'_\xi)_{ab} = \chi(ab) (\zeta^{2ab} - \zeta^{-2ab}) \hs{5ex}\nd\hs{5ex} (T'_\xi)_{ab}=\delta_{ab}\zeta^{-a^2}
\end{equation}
for $1\le a, b\le \frac{r-1}{2}$ where $\chi$ is the quadratic character of $\bbz/r$. We assume $\chi(0)=0$.

Before we compare $\phi_\zeta$ with $\phi'_\xi$, we need some basic facts from elementary number theory.

\subsection{\boldmath The field $R_r$}\label{field}
Since $r\equiv 2 (3)$, $R_r = \zo/r$ is a field of $r^2$ elements.
Recall that $u\in R_r^*$ is a generator of $\ur$ such that $u=x^{r-1}$ for a generator of $R_r^*$. Then we have
\begin{equation}\label{s1}
s := \sum_{j=0}^r(-1)^j\zeta^{\tr(u^j)} = \frac{1}{r-1}\sum_{\alpha\in R_r^*}(\zeta^{\tr(\alpha^{2r-2})} - \zeta^{\tr(\alpha^{r-1})}).
\end{equation}
Elements in $R_r$ can be written as $a-b\omega$ with $a,b\in\bbz/r$.
For any $\alpha = a-b\omega\in R^*_r$ we have
$$
\tr(\alpha^{r-1}) = \sum_{j=0}^{r-1}\binom{r-1}{j} a^{r-1-j}(-b)^j\tr(\omega^j) = \frac{3a(a+b)}{a^2+ab+b^2}-1
$$
and
$$
\tr(\alpha^{2r-2}) - \tr(\alpha^{r-1}) = \frac{-9a(a+b)b^2}{(a^2+ab+b^2)^2}.
$$
Note that for $a-b\omega \in R_r^*$, $a^2+ab+b^2\ne 0$.
Hence
\begin{eqnarray*}
s &=& \frac{1}{r-1}\sum_{a,b=1}^{r-1}(\zeta^{\tr(a-b\omega)^{2r-2}} - \zeta^{\tr(a-b\omega)^{r-1}}) \\
&=& \frac{1}{r-1}\sum_{a,b=1}^{r-1}\zeta^{\tr(a-b\omega)^{r-1}} (\zeta^{\frac{-9(1+b/a)}{(1+b/a+(b/a)^2)(1+a/b+(a/b)^2)}}-1) \\
&=& \zeta^{-1} \sum_{a=1}^{r-1}\zeta^{\frac{3(1+a)}{1+a+a^2}}(\zeta^{\frac{-9(1+a)}{(1+a+a^2)(1+1/a+1/a^2)}}-1) \\
&=& \zeta^2\sum_{a=1}^{r-1}\zeta^{\frac{-3a^2}{a^2+a+1}}(\zeta^{\frac{-9(1+a)}{(1+a+1/a)^2}}-1) 
= \zeta^2\sum_{a=1}^{r-1}\zeta^{-3(\frac{2a+1}{a^2+a+1})^2} - \zeta^2\sum_{a=1}^{r-1}\zeta^{\frac{-3}{a^2+a+1}}.
\end{eqnarray*}
Since $\bar s = -s$ if $r\equiv 1 (4)$ and $\bar s = s$ if $r\equiv 3 (4)$, we have
\begin{equation}\label{s}
s = \left\{ \begin{array}{ll}
\zeta^{-2}\sum_{a=1}^{r-1}\zeta^{\frac{3}{a^2+a+1}} -\zeta^{-2}\sum_{a=1}^{r-1}\zeta^{3(\frac{2a+1}{a^2+a+1})^2}
& \textrm{if $r\equiv 1 (4)$}, \\
\zeta^{-2}\sum_{a=1}^{r-1}\zeta^{3(\frac{2a+1}{a^2+a+1})^2} - \zeta^{-2}\sum_{a=1}^{r-1}\zeta^{\frac{3}{a^2+a+1}}
& \textrm{if $r\equiv 3 (4)$}.
\end{array} \right.
\end{equation}
For the rest of this subsection, we assume $r$ to be any odd prime.
Recall that $\zeta = \exp(\frac{2\pi\sqrt{-1}}{r})$. For any $l\in\bbz/r$, the Gauss sum is
\begin{equation}\label{Gauss}
\sum_{k=1}^{r-1} \chi(k) \zeta^{l k} = \left\{ \begin{array}{ll}
\sqrt{r} \chi(l)
& \textrm{if $r\equiv 1 (4)$}, \\
\sqrt{-r} \chi(l)
& \textrm{if $r\equiv 3 (4)$}.
\end{array} \right.
\end{equation}
It is also well known that for $a, b, c\in \bbz/r$,
\begin{equation}\label{degree2}
\sum_{i=0}^{r-1} \chi(a i^2+bi+c) = \left\{ \begin{array}{ll}
- \chi(a)
& \textrm{if $b^2-4ac\not\equiv 0 (r)$}, \\
(r-1) \chi(a)
& \textrm{otherwise}.
\end{array} \right.
\end{equation}
For any polynomial $f(x)\in \bbz/r[x]$,
\begin{equation}\label{square}
	\sum_{x=1}^{r-1}\chi(f(x^2)) - \sum_{x=1}^{r-1}\chi(f(x)) = \sum_{x=1}^{r-1}\chi(xf(x)).
\end{equation}

Besides the standard relations given in (\ref{Gauss} -\ref{square}) we will need a more involved identity that can be found in \cite{jacobsthal}:
\begin{lemma}[\cite{jacobsthal}]\label{jacobsthal}
For any $a,b\in\bbz/r$,
$$
\sum_{x=0}^{r-1}\chi(x)\chi(x^2+ax+b) = \sum_{x=0}^{r-1}\chi(x+a)\chi(x^2-4b).
$$
\end{lemma}

\subsection{\boldmath Comparing $\phi_\zeta$ to $\phi'_\xi$}\label{compare}
We assume $r\equiv 2 (3)$ again. Suppose $r\equiv 1 (4)$.
It is known that there exist constants $c_1, c_2, c_1'$ and $c_2'$ such that
\begin{equation}\label{liftS}
\tilde\phi_\zeta : \slz\to \End(V_r) \hs{3ex} \textrm{with} \hs{3ex} \tilde\phi_\zeta(S) = c_1 S_\zeta \hs{3ex} \textrm{and}
\hs{3ex} \tilde\phi_\zeta(T) = c_2 T_\zeta
\end{equation}
and
\begin{equation}\label{liftT}
\tilde\phi'_\xi : \slz\to \End(V'_r) \hs{3ex} \textrm{with} \hs{3ex} \tilde\phi'_\xi(S) = c_1' S'_\xi \hs{3ex} \textrm{and}
\hs{3ex} \tilde\phi'_\xi(T) = c_2' T'_\xi.
\end{equation}
are honest $\slz$-representations.

Suppose $r\equiv 3 (4)$.
We can also lift $\phi_\zeta$ and $\phi_\xi'$ to honest $\slz$-representations in this case as in Equations~(\ref{liftS}) and (\ref{liftT}).
But we will also take conjugation for $\phi_\xi'$, cf
Remark~\ref{rk}. So 
$$
\tilde\phi'_\xi : \slz\to \End(V'_r) \hs{3ex} \textrm{with} \hs{3ex} \tilde\phi'_\xi(S) = c_1' \bar S'_\xi \hs{3ex} \textrm{and}
\hs{3ex} \tilde\phi'_\xi(T) = c_2' \bar T'_\xi.
$$

\noindent
Our main result, Theorem~\ref{main-thm-2}, follows directly from the next theorem.

\begin{thm}\label{thm}
Suppose $r\equiv 2 (3)$.
The quantum $\psu{2}$ representation $\tilde \phi'_\xi$ is a direct summand of the
quantum $\psu{3}$ representation $\tilde \phi_\zeta$.
\end{thm}

\begin{proof}
From Equations~(\ref{eq:PSU(3)}) and (\ref{eq:PSU(2)}) we see that $T_\zeta = T'_\xi$. Lemma~\ref{S-matrix} below says that the $S$-matrices are proportional. Hence $\tilde \phi'_\xi$ is a sub-representation of $\tilde \phi_\zeta$. 
But since they are unitary representations, the theorem follows.
\end{proof}

\begin{rmk}
For any odd prime $r$, $\tilde\phi'_\xi$ is irreducible.
This can be proved using the same argument as in \cite{roberts} because the
$\psu{2}$ TQFT also has skein presentations. 
\end{rmk}

\begin{rmk}
This theorem is not true when $r\equiv 1 (3)$. For example when $r=7$, $\tilde \phi_\zeta$ is a direct sum of irreducible summands of dimensions
1 and 4 while $\tilde\phi'_\xi$ is irreducible of dimension 3.
\end{rmk}

\begin{lemma}\label{S-matrix}
$S_\zeta$ and $S'_\xi$ are proportional.
\end{lemma}

\begin{proof}
We will work out details for $r\equiv 1 (4)$ and leave the other case for the reader.
We only have to show that $S_\zeta S'_\xi = \mathrm{constant}\cdot I$.
Let $s_i := \sum_{j=0}^r (-1)^j\zeta^{i\tr(u^j)}$. Note that $s_1 = s$ as in Equation~(\ref{s1}).
Since $r\equiv 1 (4)$, $\chi(-1)=1$ and $\bar s_i = -s_i$. We have
\begin{eqnarray*}
\sum_{k=1}^{(r-1)/2}\chi(kj)s_{ik}\zeta^{-2jk} = \sum_{k=(1-r)/2}^{-1}\chi(-kj)\bar s_{ik}\zeta^{2jk}
= \sum_{k=(r+1)/2}^{r-1}-\chi(kj) s_{ik}\zeta^{2jk}.
\end{eqnarray*}
Therefore,
\begin{eqnarray*}
(S_\zeta S'_\xi)_{ij} = \sum_{k=1}^{(r-1)/2} \chi(kj)s_{ik}(\zeta^{2kj} - \zeta^{-2kj}) = \sum_{k=1}^{r-1}
\chi(kj)s_{ik}\zeta^{2kj}.
\end{eqnarray*}
Let $\alpha_{ija} = 2(j-i)+\frac{3i}{a^2+a+1}$ and $\beta_{ija}=2(j-i)+3i(\frac{2a+1}{a^2+a+1})^2$. By Equations~(\ref{s}) and
(\ref{Gauss}),
\begin{eqnarray*}
(S_\zeta S'_\xi)_{ij} &=& \chi(j) \sum_{a,k=1}^{r-1}\chi(k)(\zeta^{k\alpha_{ija}}- \zeta^{k\beta_{ija}})\\
&=& \sqrt{r}\chi(j)\sum_{a=1}^{r-1}\left[\chi(\alpha_{ija})-\chi(\beta_{ija})\right].
\end{eqnarray*}
Therefore, by Equation~(\ref{degree2}),
\begin{equation}\label{ii}
(S_\zeta S'_\xi)_{ii} =
\sqrt{r}\chi(i)\left[\sum_{a=1}^{r-1}\chi(3i)\chi(a^2+a+1)-(r-2)\chi(3i)\right]
=r\sqrt{r}.
\end{equation}
It remains to prove $(S_\zeta S'_\xi)_{ij}=0$ if $i\ne j$. In the rest of the prove we
assume that $i\ne j$.
Since $\alpha_{ij0} = \beta_{ij0}$, it is enough to prove
\begin{equation}\label{ab}
\sum_{a=0}^{r-1}\chi(\alpha_{ija}) = \sum_{a=0}^{r-1}\chi(\beta_{ija}).
\end{equation}
For $1\le i, j\le (r-1)/2$ let $b = \frac{2(j-i)}{3i}$. (Of course $b$ depends on $i$ and $j$.) Since $i\pm j\ne 0$, $b\ne 0$ or $-4/3$.
\begin{eqnarray*}
	\chi(3i)\sum_{a=0}^{r-1}\chi(\alpha_{ija}) &=& 
								 \sum_{a=0}^{r-1}\chi(b+(a^2+a+1)^{-1}) \\
								 &=&\sum_{a=0}^{r-1}\chi(a^2+a+1)
							\chi(b(a^2+a+1)+1).
\end{eqnarray*}
Since
$
\chi(a^2+a+1) = \chi(4(a^2+a+1)) = \chi((2a+1)^2+3),
$
\begin{eqnarray*}
	\chi(3i)\sum_{a=0}^{r-1}\chi(\alpha_{ija}) &=& \sum_{a=0}^{r-1}\chi((2a+1)^2+3)\chi(b((2a+1)^2+3)+4) \\
	&=& \sum_{x=0}^{r-1}\chi(x^2+3)\chi(b(x^2+3)+4).
\end{eqnarray*}
Let $c = \frac{4}{b}+3$. ($c\ne 0$ or 3.)
\begin{eqnarray*}
	\chi(3ib)\sum_{a=0}^{r-1}\chi(\alpha_{ija}) &=& \sum_{x=0}^{r-1}\chi(x^2+3)\chi(x^2+c)\\
	&=& \sum_{x=0}^{r-1}\chi(x+3)\chi(x+c) + \sum_{x=0}^{r-1}\chi(x)\chi(x+3)\chi(x+c)
\end{eqnarray*}
by Equation~(\ref{square}).
Similarly
$$
	\chi(3ib)\sum_{a=0}^{r-1}\chi(\beta_{ija}) = \sum_{x=0}^{r-1}\chi((x-3)^2+4cx) + \sum_{x=0}^{r-1}\chi(x)\chi((x-3)^2+4cx).
$$
Since $c\ne 0$ or 3,
$$
\sum_{x=0}^{r-1}\chi(x+3)\chi(x+c) = \sum_{x=0}^{r-1}\chi((x-3)^2+4cx)
$$
by Equation~(\ref{degree2}). Hence to prove Equation~(\ref{ab}), it is enough to show
\begin{equation*}
\sum_{x=0}^{r-1}\chi(x)\chi(x+3)\chi(x+c) = \sum_{x=0}^{r-1}\chi(x)\chi((x-3)^2+4cx),
\end{equation*}
which is the same as
$$
\sum_{x=0}^{r-1}\chi(x-c)\chi(x-c+3)\chi(x) = \sum_{x=0}^{r-1}\chi(x+3-2c)\chi(x^2+4c(3-c)).
$$
This equality follows from Lemma~\ref{jacobsthal} by taking $a = 3 - 2c$ and $b = c^2 - 3c$. 
This completes the proof for $r\equiv 1 (4)$.
\end{proof}

\end{document}